\documentclass[final,5p,twocolumn,times]{elsarticle}

\usepackage{graphicx,epsfig,amssymb,amsmath,amsthm,amscd,amsfonts,bm,accents} 
\usepackage[pdftex]{hyperref}
\hypersetup{
	colorlinks   = true, 
	urlcolor     = blue, 
	linkcolor    = blue, 
	citecolor   = red 
}
\usepackage{natbib}
\usepackage[labelformat=simple]{subcaption}

\usepackage{algorithm}
\usepackage{algorithmic}

\newlength\myindent
\usepackage{xcolor}
\usepackage[noabbrev]{cleveref}

\newtheorem{rmrk}{Remark}

\usepackage{mathrsfs}

\usepackage{enumitem}

\newcommand{\R}{\mathbb{R}}
\newcommand{\Pol}{\mathbb{P}}

\newcommand{\Par}{\mathscr{P}}
\newcommand{\Tk}{\mathscr{T}}

\newcommand{\Fk}{\mathscr{F}}

\usepackage{enumitem}

\newcommand{\llbrace}{\lbrace \hspace{-.045cm}\lbrace }
\newcommand{\rrbrace}{\rbrace \hspace{-.05cm}\rbrace }
\DeclareMathOperator*{\argmin}{arg\,min}

\journal{MRC}

\begin{document}
\begin{sloppypar}
\begin{frontmatter}


\title{A nonlinear weak constraint enforcement method for advection-dominated diffusion problems}

\cortext[correspondingauthor]{Corresponding author.}
\author[1,3]{Roberto J. Cier\corref{correspondingauthor}}
\ead{rcier93@gmail.com}

\author[2]{Sergio Rojas}
\ead{srojash@gmail.com}

\author[2,3]{Victor M. Calo}
\ead{vmcalo@gmail.com}
\address[1]{School of Civil and Mechanical Engineering,  Curtin University, Kent Street, Bentley, Perth, WA 6102, Australia}
\address[2]{School of Earth and Planetary Sciences, Curtin University, Kent Street, Bentley, Perth, WA 6102, Australia}
\address[3]{Mineral Resources, Commonwealth Scientific and Industrial Research Organisation (CSIRO), Kensington, Perth, WA 6152, Australia}




\begin{abstract}
  
  We devise a stabilized method to weakly enforce bound constraints in the discrete solution of advection-dominated diffusion problems.  This method combines a nonlinear penalty formulation with a discontinuous Galerkin-based residual minimization method. We illustrate the efficiency of this scheme for both uniform and adaptive meshes through proper numerical examples.

\end{abstract}


\end{frontmatter}


\section{Introduction}
Standard (Galerkin) finite element methods (FEM) can yield unphysical oscillatory discrete solutions in advection-dominated regimes. A commonly used technique to overcome this weakness of the formulation is to add stabilized terms that enhance the properties of the discrete solution. Some of these techniques yield Petrov-Galerkin schemes, such as the SUPG method~\cite{ brooks1982} or the streamline diffusion (SD) method~\cite{ johnson1986}. Other stabilization techniques include least-squares formulations~\cite{ bramble1997}, variational multiscale (VMS)~\cite{ hughes2007}, subgrid viscosity~\cite{ guermond1999}, and continuous interior penalty methods (CIP)~\cite{ burman2004}, among others. 

Although stabilized formulations improve the robustness and accuracy of the numerical solutions, spurious undershoots and overshoots are not eliminated, especially in low-resolution meshes. These oscillations are a drawback as in many engineering applications (i.e., transport of density, concentration, or temperature) require to remain within their physical range. Violating these bounds delivers poor simulation outputs. Thus, overshoots or undershoots are controlled through proper constraint enforcement procedures. For that reason, a plethora of techniques to surmount this effect has been proposed, mostly constructed from a stabilized formulation. One of these schemes incorporates shock-capturing terms to satisfy a discrete maximum principle~\cite{ burman2005, mizukami1985}. Also, flux-corrected methods~\cite{kuzmin2002, lohner1987} seek to impose the constraints by altering the system matrix. These methods are generally only first-order accurate. Higher-order schemes require terms to control and, in many times, reduce the dissipative response of the method. 

More recently, an alternative constraint imposition approach --more precisely, positivity preserving-- was proposed in~\cite{ burman2017}. The authors satisfy the discrete maximum principle weakly by adding a consistent penalty term to the variational formulation of a Galerkin least-squares (Ga-LS) finite element discretization. This method is flexible and can incorporate a priori lower and upper bounds on the discrete solution, by simply adding the corresponding consistent penalty term to the discrete formulation. We combine this consistent penalization with a new adaptive stabilized finite element framework that minimizes the residual in dual norms of discontinuous Galerkin~(dG) methods~\cite{calo2019}. This formulation inherits the stability and accuracy of the underlying dG approximation. The formulation seeks for a solution in a continuous trial function space which is a proper subspace of the dG function space. The resulting saddle-point problem delivers stable formulations with continuous solutions with a robust \emph{a posteriori} error estimate, which can be computed on the fly to drive optimal adaptive mesh refinements.

In this paper, we develop a constraint enforcement technique that combines the ideas of the nonlinear penalty method of~\cite{ burman2017}, and the residual minimization technique of~\cite{ calo2019}, applied to advection-dominated diffusion problems. We construct it as follows. First, we modify the corresponding bilinear dG form by adding a nonlinear penalty term to enforce constraints weakly. Next, we solve a residual minimization problem in a dG dual norm. The resulting technique minimizes weakly the violation of solution bounds and additionally delivers a robust residual estimator to guide adaptive mesh refinement. The main advantage of considering this procedure is that it results in a nonlinear saddle-point problem, with symmetric Jacobian. Therefore, an extensive list of iterative solvers is available for each step of the Newton iteration (see, e.g., \cite{benzi2005}). The idea of combining residual minimization, together with nonlinear techniques, was also considered in~\cite{muga2019}, as an extension to advection-reaction problems in Banach spaces, and in~\cite{houston2020} as a technique to remove the Gibbs phenomena in diffusion-advection-reaction problems. However, the main difference with this work is that, in those, the nonlinearity appears in the dual norm.

The paper goes as follows. In Section~\ref{section:model}, we state the model problem and the dG discretization. We introduce the consistent penalty method in Section~\ref{section:penalty} and we derive the nonlinear formulation for weak constraint enforcement. Besides, we detail the adopted resolution scheme for the linearization of the problem. Finally, in Section~\ref{section:numexp}, we report some numerical experiments illustrating the solution quality with mild and sharp inner layers, as well as discontinuities, with diminished violations of the discrete solution bounds.

\section{Model problem: Diffusion--advection-reaction}\label{section:model}

In this section, we present all the required ingredients for the constraint enforcement method in the context of advection-diffusion-reaction problems. Let $\Omega$ be an open, bounded Lipschitz domain in ${\R^d, d\in \{2, 3\}}$, with boundary ${\partial\Omega}$; ${\beta: \Omega \rightarrow \R^d}$ be an advective velocity field; ${K \in \R^{d\times d}}$ be a diffusion tensor, assumed to be continuous, symmetric and positive-definite; and ${\sigma: \Omega \rightarrow \R}$ be a reaction coefficient. We assume that ${f: \Omega \rightarrow \R}$ is a given source term, and ${g: \Gamma \rightarrow \R}$ is a prescribed Dirichlet boundary condition. We consider the following advection-diffusion-reaction problem:
\begin{equation}\label{eq:adr_problem}
\hspace{-0.225cm}\left\{\begin{array}{l}
\text{Find } u \text{ such that:} \smallskip \\
\begin{array}{rcll}
\hspace{-0.2cm}A(u):=-\nabla \cdot \left( K \nabla u \right) + \beta \cdot \nabla u + \sigma \, u &\hspace{-0.25cm}=& \hspace{-0.25cm}f,& \hspace{-0.35cm}\text{ in } \Omega, \smallskip\\
u &\hspace{-0.25cm}=&\hspace{-0.25cm} g, & \hspace{-0.35cm}\text{ on } \Gamma,\smallskip
\end{array}
\end{array}\right.
\end{equation}
where, using the standard notation, we assume that ${\beta\in [W^{1, \infty}(\Omega)]^d}$, ${\sigma \in L^{\infty}(\Omega)}$, ${f\in L^2(\Omega)}$ and ${g \in H^{1/2}(\partial\Omega)}$. Then, the weak formulation of~\eqref{eq:adr_problem} becomes:
\begin{equation}\label{eq:weak_form}
\hspace{-0.25cm}\left\{\begin{array}{l}
\text{Find } u \in H^1(\Omega), \text{ such that:} \smallskip \\
(K\nabla u, \nabla v)_{0, \Omega} + (\beta \cdot \nabla u, v)_{0, \Omega} + (\sigma u, v)_{0, \Omega} = (f, v)_{0, \Omega}, \\
\forall v \in H^1_0(\Omega)
\end{array}\right.
\end{equation}
where $(\cdot,\cdot)_{0,\Omega}$ denotes the $L^2$-scalar product on $\Omega$. Henceforth, we assume that there is a real number ${\sigma_0 \geq 0}$, such that ${\sigma - \frac{1}{2}\nabla \cdot \beta \geq \sigma_0}$ in $\Omega$. Owing to the above assumptions, the Lax-Milgram Lemma implies that problem~\eqref{eq:weak_form} is well posed~\cite{ ciarlet2002}. In what follows, we assume that the exact solution $u$ is in $H^2(\Omega)$. If the reaction source satisfies ${\sigma\geq0}$, problem~\eqref{eq:weak_form} also satisfies a maximum principle, that is, under suitable assumptions on the data $f$ and $g$, the solution attains its maximum or minimum at the boundary. In particular, if ${f\geq0}$ and ${g\geq0}$, then ${u(x)\geq0, \forall x \in \Omega}$. Similarly, in pure convection-diffusion problems (i.e., if ${f=0}$ and ${\sigma=0}$), then ${ \min_{y\in \partial\Omega}g(y)\leq u(x)\leq \max_{y\in \partial\Omega}g(y), \forall x \in \Omega}$. For a detailed discussion on maximum principles for elliptic second-order problems see~\cite{ gilbarg2015}.

Given that in this work we focus on advection-dominated cases, that is, when ${\|\sigma\|_\infty, \|K\|_\infty \ll \|\beta\|_{\ell}}$, where $\ell$ is a length scale, we conveniently split the boundary ${\partial\Omega\equiv\Gamma}$ into ${\Gamma=\Gamma_-\cup \Gamma_0\cup \Gamma_+}$, with
\begin{equation*}
\begin{array}{ll}
\Gamma_-&= \ \{x\in{\Gamma} \ ; \ \beta\cdot n<0\} \quad \textrm{(inflow boundary)},\\
\Gamma_0&= \ \{x\in{\Gamma} \ ; \ \beta\cdot n=0\} \quad \textrm{(characteristic boundary)}, \\
\Gamma_+&= \  \{x\in{\Gamma} \ ; \ \beta\cdot n>0\} \quad \textrm{(outflow boundary)},
\end{array}
\end{equation*}
where $n$ represents the outward normal vector to $\Gamma$.

\subsection{Discontinuous Galerkin variational formulation}\label{section:dG}

In this section, we describe the dG variational formulation that we enlarge by including the constraint enforcement penalty terms. 

Let $\{\Tk_h\}$ be a family of simplicial meshes of $\Omega$. For simplicity, we assume that any mesh exactly represents $ \Omega$ in $\Tk_h$, that is, $\Omega$ is a polygon or a polyhedron. $ T$ denotes a generic element in $\mathscr{T}_h$, $h_T$ denotes the diameter of $T$ and $n_T$ its outward unit normal. We set $h = \max_{T \in \mathscr{T}_h} h_T$. We assume, without loss of generality, that $h \leq 1$. We define the classical dG approximation space
\begin{equation*}
V_h :=\{v_h \in L^2(\Omega) \ | \  \forall T \in \mathscr{T}_h, v_h|_T \in \Pol^p \},
\end{equation*}
where $\Pol^p$ denotes the set of polynomials, defined over $T$, with polynomial degree smaller or equal than $p$. 

Let $F$ be an interior face of the mesh if there are $T^-(F)$ and $T^+(F)$ in $\mathscr{T}_h$, such that ${F=T^+(F) \cap T^-(F)}$, and we let $n_F$ be the unit normal vector to $F$ pointing from $T^-(F)$ towards $T^+(F)$. Similarly, $F$ is a boundary face if there is a ${T(F)\in \mathscr{T}_h}$ such that ${F=T(F)\cap \Gamma}$, and we let $n_F$ coincide with $n$. We collect all the faces or edges into the set ${\mathscr{F}_h= \bigcup_{T\in\mathscr{T}_h} F}$. We define the boundary skeleton $\mathscr{F}^{b}_h$ as ${\mathscr{F}^{b}_h = \mathscr{F}_h \cap \Gamma}$, and the internal skeleton $\mathscr{F}^i_h$ as ${\mathscr{F}^i_h=\mathscr{F}_h \backslash \Gamma}$. Henceforth, we deal with functions that are double-valued on $\mathscr{F}^i_h$ and single-valued on $\mathscr{F}^{b}_h$, for example, all functions in $V_h$ have these characteristics. On interior faces, when the two branches of the function in question, say $v$, are associated with restrictions to the neighboring elements $T^\mp(F)$, we denote these branches by $v^\mp$, and the jump $[\![v]\!]_F$ and the standard (arithmetic) average $\{\!\{v\}\!\}_F$ as 
\begin{align*}
[\![v]\!]_F &:= v^- - v^+, \quad \quad \{\!\{v\}\!\}_F := \frac{1}{2}(v^- + v^+),
\end{align*}
On a boundary face $F \in \mathscr{F}^b_h$, we set $[\![v]\!]_F=\{\!\{v\}\!\}_F=v|_F$. The subscript $F$ is omitted from the jump and average operators when there is no ambiguity. Finally, we set $h_F$ as the diameter of the face $F$.

Given the previous components, the dG discretized formulation for~\eqref{eq:adr_problem} reads:
\begin{equation}\label{eq:gen_vf}
\left\{\begin{array}{l}
\text{Find } u_h \in V_h, \text{ such that:} \smallskip \\
\begin{array}{l}
b_h(u^{\text{dG}}_h,v_h)= \ell_h(v_h), \quad \forall v_h \in V_h,
\end{array}
\end{array}\right.
\end{equation}
with the bilinear form
\begin{align*}
b_h(u^{\text{dG}}_h,v_h):=b_h^\textrm{diff}(u^{\text{dG}}_h,v_h) + b_h^\textrm{adv}(u^{\text{dG}}_h,v_h),
\end{align*}
where
\begin{equation*}
\begin{array}{rcl}
\hspace{-1.cm} \displaystyle b_h^\textrm{diff}(w,v) &:= & \displaystyle \sum_{T \in \Tk_h} (K \nabla w \, , \, \nabla v)_{0,T} \\
\hspace{-1.cm}& +  &  \displaystyle \sum_{F \in \Fk_h}  \theta \, \left(  \llbracket w \rrbracket \, , \,  \llbrace  K \nabla v \rrbrace \cdot n_F \right)_{0,F} , \smallskip\\
\hspace{-1.cm}& -  &  \displaystyle  \sum_{F \in \Fk_h}\left(  \llbrace  K \nabla w \rrbrace \cdot n_F \, , \,  \llbracket v \rrbracket \right)_{0,F} \\
\hspace{-1.cm}& +  & \displaystyle \sum_{F \in \Fk_h} \eta \left( \llbracket w \rrbracket \, , \, \llbracket v \rrbracket \right)_{0,F} \\
\end{array}
\end{equation*}
\begin{equation*}
\begin{array}{rcl}
\hspace{-1.cm} \displaystyle b_h^\textrm{adv}(w,v) & := & \displaystyle \sum_{T\in\Par_h}(\, \beta \cdot \nabla  w + \sigma \,  w\, , \, v)_{0,T} \smallskip                                                                                                                                                      \\
\hspace{-1.cm}& +  &  \displaystyle \sum_{F\in\Fk^b_h\cap \, \Gamma^-} \left( \, \beta \cdot n_F \, w, v\right)_{0,F} \smallskip \\
\hspace{-1.cm}& -  &  \displaystyle \sum_{F \in \Fk^i_h} \left( \, \beta \cdot n_F \, \llbracket w \rrbracket \, , \, \llbrace v \rrbrace \right)_{0,F} \smallskip \\
\hspace{-1.cm}& +  &   
\displaystyle \sum_{F \in \Fk^i_h}\dfrac{1}{2} \left(\, \left| \, \beta \cdot n_F\right| \, \llbracket w \rrbracket \, , \, \llbracket v \rrbracket \right)_{0,F}, \smallskip
\end{array}
\end{equation*}
and the linear form
\begin{equation*}
\begin{array}{rcl}
\hspace{-0.2cm} \displaystyle \ell_h(v) &              \hspace{-0.2cm}:=               & \hspace{-0.2cm}\displaystyle \sum_{T\in\Tk_h}(f, v)_{0,T} +  \sum_{F \in \Fk_h^b} \eta \left( g, v \right)_{0,F}  \\
&             \hspace{-0.2cm} +                 &\hspace{-0.2cm}\displaystyle \sum_{F\in\Fk^b_h\cap \, \Gamma^-} \left( \, \beta \cdot n_F \, g,  v\right)_{0,F} + \sum_{F\in\Fk^b_h} \theta \, \left( g,  K \nabla v \cdot n \right)_{0,F}. \smallskip
\end{array}
\end{equation*}

For diffusion problems, we recover well-known types of dG formulations for different choices of $\theta$  and the penalty parameter $\eta$ in $b^{\textrm{diff}}_h(w, v)$ (e.g., see~\cite{arnold2002unified, riviere2008}). Herein, for our numerical experiments we set the parameters to deliver the SIPG method, that is, ${\theta=-1}$ and, following~\cite{ shahbazi2005}, $\eta = \eta_0(p+1)(p+d)K/h$, being $p$ the polynomial degree for the test functions and $\eta_0=3$. In the advective part of the bilinear form, $b^{\text{adv}}_h(w, v)$, we use an upwinding scheme (see~\cite{ brezzi2004, ern2011}). The broken polynomial space $V_h$ can be equipped with the following norm:
\begin{align*}\label{dg_norm}
\|w\|^2_{V_h} := \|w\|^2_{\text{adv}} + \|w\|^2_{\text{diff}},
\end{align*}
with $\|w\|^2_{\text{adv}}$ representing the upwinding norm defined for advection-reaction problems and $\|w\|^2_{\text{diff}}$ representing the norm defined by the interior penalty methods for diffusion problems. Thus, these norms read
\begin{equation*}
\begin{array}{rcl}
\hspace{-0.2cm} \|w\|^2_{\text{adv}} &\hspace{-0.2cm}:=&\hspace{-0.2cm} \displaystyle \|w\|^2_{0, \Omega} + \frac{1}{2}\| \, |\, \beta\cdot n|^{\frac{1}{2}}w\|^2_{0, \Gamma} \\
&\hspace{-0.2cm}+&\hspace{-0.2cm} \displaystyle \frac{1}{2}\sum_{ F \in \mathscr{F}^i_h} \left(\, | \, \beta\cdot n_F| \, [\![w]\!], [\![w]\!] \right) _{0,F} \\
&\hspace{-0.2cm}+&\hspace{-0.2cm} \displaystyle \sum_{T \in \mathscr{T}_h} h_T \| \, \beta\cdot\nabla{w}\|^2_{0, T}, \\
\hspace{-0.2cm}  \|w\|^2_{\text{diff}} &\hspace{-0.2cm}:=&\hspace{-0.2cm} \displaystyle \|K^{\frac{1}{2}}\nabla_h w\|^2_{0, \Omega} +  \sum_{F \in \mathscr{F}_h} \left(\eta[\![w]\!], [\![w]\!] \right) _{0,F},
\end{array}
\end{equation*}

\section{Weak constraint enforcement method based on residual minimization}\label{section:penalty}

For the sake of simplicity, in the next we assume that the aim is to enforce a positivity preserving condition, that is, ${u \geq 0}$. Other varieties of constraints, such as upper bounds or other minimal values, can also be imposed by considering a slight modification of the nonlinear term (see Remark~\ref{upper_lower_bounds}). 

\subsection{Nonlinear consistent penalty method}
Consider the following penalization term (see~Remark~\ref{rmrk2}):
\begin{equation}\label{eq:penalty}
	\displaystyle {\gamma=\gamma_0 \left(\frac{\|\,\beta\,\|_\ell}{h}+\frac{\|K\|_\infty}{h^2} + \|\sigma\|_\infty \right)^{-1}},
	\end{equation}
where $0<\gamma_0<1$ is a user-defined constant real number.
We define ${\xi_\gamma: V_h \rightarrow \R}$, as the function:
\begin{align}
\xi_\gamma(v_h):=[v_h-\gamma(A(v_h)-f)]_-, \quad \forall v_h \in V_h,
\end{align} 
where  ${x_-=\frac{1}{2}(x-|x|)}$ denotes the negative part of the real number $x$, satisfying  ${x_{-}=x}$ if ${x<0}$, and ${x_- = 0}$ if ${x\geq 0}$. 

We define ${b^\gamma_h(u_h;v_h)}$, composed by the original bilinear form ${b_h(u_h,v_h)}$ and a nonlinear penalty term, as follows:
\begin{align}\label{nonlinear_form}
b^\gamma_h(u_h;v_h) := b_h(u_h, v_h) + \langle \gamma^{-1} \xi_\gamma(u_h),v_h\rangle_h,
\end{align}
where
\begin{align*}
\langle x_h,y_h\rangle_h:= \sum_{ T \in \mathscr{T}_h}(x_h,y_h)_T.
\end{align*}
%
By construction, the analytical solution satisfies that $\xi_\gamma(u)=0$, since $A(u)=f$ and $u_-=0$. We consider the following discrete problem: 
\begin{equation}\label{eq:dg+penalty}
\left\{\begin{array}{l}
\text{Find } u_h \in V_h, \text{ such that:} \smallskip \\
\begin{array}{l}
b^\gamma_h(u_h;v_h)= \ell_h(v_h), \quad \forall v_h \in V_h.
\end{array}
\end{array}\right.
\end{equation}

Since $\xi_\gamma(u)$ vanishes identically in $\Omega$, exact consistency holds for~\eqref{eq:dg+penalty}. Consistency still holds if we substitute the penalty parameter $\gamma$ by a function taking uniformly positive values in $\Omega$.
\begin{rmrk}\label{upper_lower_bounds}
The nonlinear form $b^\gamma_h(u_h; v_h)$ can also impose a constraint on the upper limit of the solution. For instance, if it is known that $u\in [u_\emph{min},u_\emph{max}]$, $b^\gamma_h(u_h;v_h)$ can be written as:
\begin{equation}\label{b_gamma_01}
\begin{array}{rcl}
b^\gamma_h(u_h; v_h)&:=&\displaystyle b_h(u_h, v_h) + \langle \gamma^{-1}\xi^\emph{min}_\gamma(u_h), v_h\rangle_h \\
&+& \displaystyle \langle \gamma^{-1}\xi^\emph{max}_\gamma(u_h), v_h \rangle_h,
\end{array}
\end{equation}
 with ${\xi^\emph{min}_\gamma(u_h):=[(u_h-u_\emph{min})-\gamma(A(u_h)-f)]_-}$ and ${\xi^\emph{max}_\gamma(u_h):=[(u_\emph{max}-u_h)-\gamma(A(u_h)-f)]_-}$ representing the penalty terms imposed for controlling the lower and upper limits of the solution, respectively.
\end{rmrk}
\begin{rmrk}\label{rmrk2}
The election of the stabilization term \eqref{eq:penalty} is motivated by the classical stabilization parameters (SUPG, Ga-LS, VMS) for diffusive problems (see~\cite{codina2000}), and the stabilization parameter considered in~\cite{burman2017} for advective problems. Naive elections of the stabilization term, such as $\gamma$ constant, could affect the convergence of the solution. 
\end{rmrk}

\subsection{Discontinuous Galerkin-based residual minimization method}

We apply the adaptive stabilized method introduced in~\cite{calo2019} to diffusion-advection-reaction problems. We seek the discrete solution in a continuous trial space (e.g., $H^1$-conforming finite elements) as the minimizer of the residual measured in a suitable dG space. This procedure inherits all the desirable stability properties from the well-posed dG variational formulation. In practice, such a residual minimization leads to a stable saddle-point problem involving the continuous trial space and the discontinuous test space. The discrete solution delivers a residual representative that is an efficient and reliable error estimate to drive adaptive mesh refinement. Thus, we compute on the fly a discrete solution in the continuous trial space and an error representation in the discontinuous test space.

Starting from the stable dG formulation of the form~\eqref{eq:gen_vf} in $V_h$, a trial subspace $U_h \subset V_h$ is chosen and, rather than solving the typical square problem in $V_h$, we state the following residual minimization: 
\begin{equation}\label{eq:min_prob}
\left\{\begin{array}{l}
\text{Find } u_h \in U_h \subset V_h,  \text{ such that:} \smallskip \\
\begin{array}{rcl}
\displaystyle u_h &=& \displaystyle \argmin_{z_h \in U_h} \dfrac{1}{2}\|\ell_h- B_h \, z_h\|^2_{V_h^\ast} \\
&=& \displaystyle \argmin_{z_h \in U_h} \dfrac{1}{2}\|R^{-1}_{V_h}(\ell_h- B_h z_h)\|^2_{V_h},
\end{array}
\end{array}
\right.
\end{equation}
 where the dual norm $\| \cdot \|_{V^*_h}$ for $\varphi \in V^*_h$ is:
\begin{align}
\| \varphi \|_{V^*_h} := \sup_{v_h \in V_h \backslash \{0\}} \frac{\langle \varphi, v_h \rangle_{V^*_h \times V_h}}{\| v_h \|_{V_h}},
\end{align}
 and ${B_h: U_h \rightarrow V_h^*}$ is:
\begin{align}\label{B_operator}
\langle B_h z_h, v_h\rangle_{V_h^* \times V_h} := b_h(z_h, v_h),
\end{align} 
$\langle \cdot, \cdot \rangle_{V_h^* \times V_h}$ denotes the duality pairing in $V_h^* \times V_h$, and $R_{V_h}^{-1}$ denotes the inverse of the Riesz map:
\begin{align}\label{Riesz}
R_{V_h} \quad : & \quad V_h \rightarrow V_h^*& & \nonumber \\
& \langle R_{V_h}y_h, v_h \rangle_{V_h^* \times V_h}:=(y_h, v_h)_{V_h}. &
\end{align}

The second equality in~\eqref{eq:min_prob} holds, since the Riesz map is an isometric isomorphism. With all the above, it can be shown that~\eqref{eq:min_prob} is equivalent to the following saddle-point problem (see~\cite{ calo2019}):
\begin{equation}\label{eq:saddle-point}
\hspace{-0.5cm} \left\{\begin{array}{l}
\text{Find } (\varepsilon_h, u_h) \in V_h \times U_h,  \text{ such that:} \smallskip \\
\begin{array}{rcll}
(\varepsilon_h,v_h)_{V_h} + b_h(u_h,v_h) \hspace{-0.2cm}&=\hspace{-0.2cm}& \ell_h(v_h),& \hspace{-0.2cm}\forall v_h \in V_h, \smallskip\\
b_h(z_h, \varepsilon_h) \hspace{-0.2cm}& = \hspace{-0.2cm}&0,& \hspace{-0.2cm}\forall z_h \in U_h,
\end{array}
\end{array}\right.
\end{equation}

According to~\cite{calo2019}, the well-posedness of the dG-based residual minimization method relies on the classical assumptions for well-posedness of the original dG formulation (i.e., consistency, boundedness, and stability). The residual representative is efficient and reliable under a suitable saturation assumption. The resulting linear system leads to a saddle-point problem irrespective of the symmetry properties of the dG variational formulation, opening the possibility to use efficient well-known iterative solvers for its resolution.

\subsection{Extension for the nonlinear penalty method}

In this section, we extend the discrete formulation to solve a nonlinear problem of the form: ${N_h(u_h)=\ell_h}$, where ${N_h : U_h \rightarrow V_h^*}$ represents the operator that includes the nonlinear penalty term, defined as ${\langle N_h(z_h), v_h\rangle_{V_h^* \times V_h} := b^\gamma_h(z_h; v_h)}$. Given that ${b^\gamma_h(z_h; v_h)}$ is built from the original bilinear form, the discrete problem~\eqref{eq:dg+penalty} presents unique solution.

At the discrete level, we seek a minimizer $u_h \in U_h \subset V_h$ for the residual $\ell_h - N_h (z_h) $ associated to~\eqref{eq:dg+penalty}:
\begin{equation}\label{eq:resmin_nonlin}
\left\{\begin{array}{l}
\text{Find } u_h \in U_h \subset V_h,  \text{ such that:} \smallskip \\
\begin{array}{rcl}
\displaystyle u_h &=& \displaystyle \argmin_{z_h \in U_h} \dfrac{1}{2}\|\ell_h- N_h \, (z_h) \|^2_{V_h^\ast} \\
&=& \displaystyle \argmin_{z_h \in U_h} \dfrac{1}{2}\|R^{-1}_{V_h}(\ell_h- N_h (z_h))\|^2_{V_h},
\end{array}
\end{array}
\right.
\end{equation}

Similar to~\eqref{eq:min_prob}, we state the nonlinear problem as a critical point of the minimizing functional, which translates into the following linear problem:
\begin{equation}\label{eq:resmin_nonlin_2}
\left\{\begin{array}{l}
\text{Find } u_h \in U_h \subset V_h,  \text{ such that:} \smallskip \\
(R_{V_h}^{-1}(\ell_h - N_h(u_h)),R_{V_h}^{-1}DN_h(u_h; z_h)) = 0, \forall z_h \in U_h.
\end{array}
\right.
\end{equation}
$DN_h : U_h \rightarrow V_h^{*}$ is defined as:
\begin{align}
	\langle DN_h(u_h; z_h), v_h\rangle_{V_h^* \times V_h} := db^\gamma_h(u_h; z_h, v_h),
\end{align}
 where $db^\gamma_h(u_h; z_h, v_h)$ represents the derivative of the nonlinear form $b^\gamma_h(u_h;v_h)$ in the direction of an increment $z_h$:
\begin{align}
	db^\gamma_h(u_h; z_h, v_h):=\frac{d}{d\epsilon}b^\gamma_h(u_h+\epsilon z_h; v_h)\big|_{\epsilon=0},
\end{align}
for instance, if we can impose a positivity preserving condition through the penalty term, the derivative reads:
\begin{align}
db^\gamma_h(u_h; z_h, v_h):=b_h(z_h, v_h) + \langle \frac{1}{\gamma}d\xi_{\gamma}(u_h;z_h), v_h\rangle_h
\end{align}
where ${d\xi_{\gamma}(u_h;z_h) = \frac{1}{2}[1-\text{sgn}(u_h-\gamma(Au_h-f))][z_h-\gamma Az_h]}$. Hence, the modified discrete formulation reads: 
\begin{equation}\label{eq:saddle-point-nl}
\hspace{-0.5cm} \left\{\begin{array}{l}
\text{Find } (\varepsilon_h, u_h) \in V_h \times U_h,  \text{ such that:} \smallskip \\
\begin{array}{rcll}
(\varepsilon_h,v_h)_{V_h} + b^\gamma_h(u_h;v_h) \hspace{-0.2cm}&=\hspace{-0.2cm}& \ell_h(v_h),& \hspace{-0.2cm}\forall v_h \in V_h, \smallskip\\
db^\gamma_h(u_h; z_h, \varepsilon_h)  \hspace{-0.2cm}& = \hspace{-0.2cm}&0,& \hspace{-0.2cm}\forall z_h \in U_h,
\end{array}
\end{array}\right.
\end{equation}
The first line of the system~\eqref{eq:saddle-point-nl} represents the nonlinear problem to solve, whereas the second line linearizes the constraint we seek to impose. 
\begin{rmrk}\label{rem:convergence}
In practice, solving~\eqref{eq:saddle-point-nl} implies that a price in the energy norm may be paid to enforce the constraints, since the residual minimization method without penalty achieves the lowest possible variational residual for the linear problem (see~\cite{calo2019}, Theorem 2).
\end{rmrk}

\subsection{Solution scheme}

Given the discrete solution pair ${(\varepsilon^k_h, u^k_h)}$ in an iterative step $k$, we seek for the increment ${(\delta \varepsilon_h, \delta u_h)}$ in the next iteration step, such that ${u}^{k+1}_h={u}^{k}_h+t^k{\delta u}_h$, and ${\varepsilon}^{k+1}_h={\varepsilon}^{k}_h+t^k{\delta \varepsilon}_h$, being $t^k$ a relaxation parameter. The method looks for a solution pair $(\varepsilon_h^{k+1}, u_h^{k+1})$ that accomplishes~\eqref{eq:saddle-point-nl} to first order. We propose a solution strategy that applies Newton's method to the nonlinear problem. Considering this,~\eqref{eq:saddle-point-nl} we solve the following linearized problem at the iteration $k+1$:
\begin{equation}\label{eq:nl-newton}
\left\{\begin{array}{l}
	\text{Given the pair } (\varepsilon^k_h, u^k_h), \text{ find } (\delta \varepsilon_h, \delta u_h) \in V_h \times U_h, 
	\\\text{such that:}\\ 
	\begin{array}{llll}
		\hspace{-0.25cm}(\delta \varepsilon_h, v_h)_{V_h} & \hspace{-1cm}+ \, db^\gamma_h(u_h^{k}; \delta u_h, v_h) =                            &                                      &  \\
		                                                        & \hspace{-0.25cm} \ell_h(v_h)-b^\gamma_h(u_h^{k}; v_h) - (\varepsilon^k_h, v_h)_{V_h}, & \hspace{-0.25cm}\forall v_h\in V_h,  &  \\
		\hspace{-0.25cm}db^\gamma_h(u^k_h; z_h, \delta \varepsilon_h) & \hspace{-0.25cm}=-db^\gamma_h(u^k_h; z_h, \varepsilon^k_h),                             & \hspace{-0.25cm}\forall z_h \in U_h. &
	\end{array}
\end{array}\right.
\end{equation}

The matrix formulation of~\eqref{eq:nl-newton} reads
\begin{align}\label{eq:femwdg_mtx}
\begin{pmatrix}
\ G & B_u \ \\ \ B_u^T & 0 \ 
\end{pmatrix}
\begin{pmatrix}
\ \delta \varepsilon_h \ \\
\ \delta u_h \
\end{pmatrix}
=
\begin{pmatrix}
\ L  \ \\
\ 0\
\end{pmatrix} 
-
\begin{pmatrix}
\ G\varepsilon_h^k + N(u_h^k)\ \\
\ B^T_u \varepsilon_h^k \
\end{pmatrix} 
\end{align}
where the superindex $T$ denotes transpose. $G$ represents the Gram matrix associated with the inner product which induces the norm in the discrete space $V_h$, $N(u_h^k)$ is related to the nonlinear form $b_h^\gamma(u_h;v_h)$, and $B_u$ represents the matrix associated with its linearization $db^\gamma_h(u_h^{k}; \delta u_h, v_h)$. The residual representative $\varepsilon_h$ is a function of $u_h$. We define ${\boldsymbol{ x_h}=( \varepsilon_h,  u_h)}$, comprising both the solution and the residual representative, being valid also for the increments, which allows us to rewrite~\eqref{eq:femwdg_mtx} as:
\begin{equation*}
\boldsymbol{J}^k \, \boldsymbol{\delta x_h} = \boldsymbol{R}^k,
\end{equation*}
where 
\begin{equation*}
\begin{array}{c}
\boldsymbol{J}^k =
\begin{pmatrix}
\ G & B_u \ \\ \ B_u^T & 0 \ 
\end{pmatrix}
\, \; \text{ and } \; \,
\boldsymbol{R}^k
=
\begin{pmatrix}
\ L  \ \\
\ 0\
\end{pmatrix} 
-
\begin{pmatrix}
\ G\varepsilon_h^k + N(u_h^k)\ \\
\ B^T_u \varepsilon_h^k \
\end{pmatrix} 
\end{array}
\end{equation*}
The convergence of the method depends on the step size. Thus, we use a damped Newton algorithm to control convergence~\cite{bank1981}, with ${\omega=0.5}$ as damping parameter (see Algorithm~\ref{alg:1}). Presently, we cannot provide a bound on the number of iterations the proposed algorithm needs to achieve convergence. Nevertheless, in our experience, the algorithm is efficient and has a reasonable cost compared to the original linear problem, as shown in the following Section.
\begin{algorithm}
	\caption{Damped Newton method.}
	\label{alg:1}
	(1) input $\omega \in (0,1)$, $\zeta=0$, $k=0$, $\text{TOL}$\\
	(2) compute $\boldsymbol{x_h^k}=(\varepsilon_h^k, u_h^k)$\\
	(3) compute $\|\boldsymbol{R}^k\|$ \\
	(4) $\displaystyle t^k = \frac{1}{1 + \zeta \|\boldsymbol{R}^k\|}$ \\
	(5) compute $\displaystyle \boldsymbol{x_h}^{k+1} = \boldsymbol{x_h} + t^k\boldsymbol{\delta x_h}$, $\|\boldsymbol{R}^{k+1}\|$\\
	(6) \textbf{if} $\displaystyle \frac{1}{t^k}\left(1-\frac{\|\boldsymbol{R}^{k+1}\|}{\|\boldsymbol{R}^{k}\|}\right) < \omega$\\
	(7) \textbf{then} \{\textbf{if} $\zeta=0$ \textbf{then} $\zeta=1$ 
	\textbf{else} \{$\zeta=10\zeta$; go to (4)\}\} \\
	(8) \textbf{else} \{$\zeta=\zeta/10$; $k=k+1$\}\\
	(9) \textbf{if} $\|u^{k+1}-u^k\|<\text{TOL}$ \textbf{then} return \textbf{else} go to (3)	
\end{algorithm}

\section{Numerical experiments}\label{section:numexp}

In this section, we implement the nonlinear constraint enforcement method to solve several numerical tests using FEniCS~\cite{ alnaes2015}. 

\subsection{Advection problem over a quasi-uniform mesh}

\begin{figure}[h!]
  \centering
  \begin{subfigure}{0.35\textwidth}
    \centering
    \includegraphics[trim= 3cm 1.5cm 0cm 2cm,width=\textwidth]{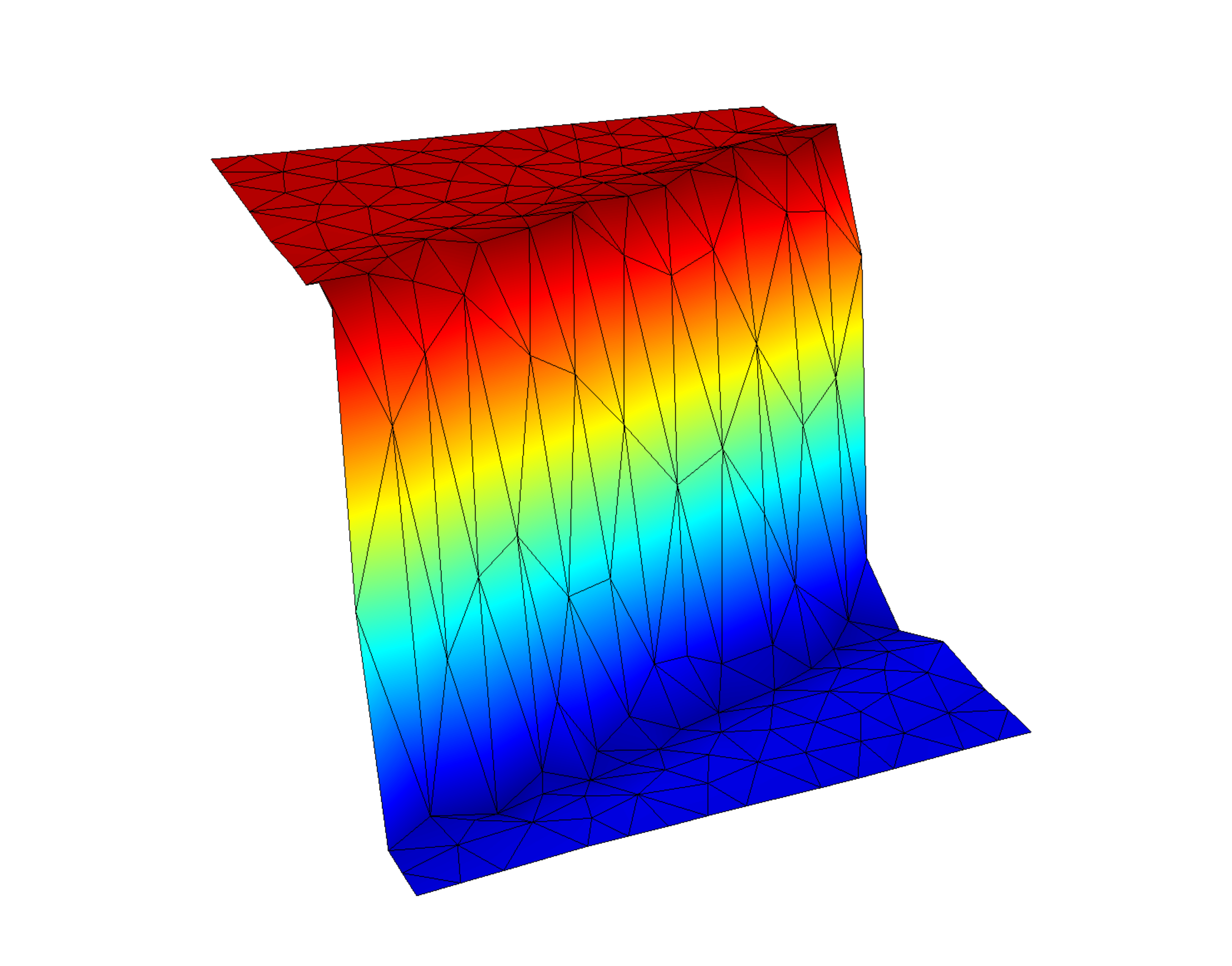} 
    \caption{\footnotesize Solution without penalty method.}
    \label{fig:tc1-sin}
  \end{subfigure}
  \begin{subfigure}{0.35\textwidth}
    \centering
    \includegraphics[trim= 3cm 1.5cm 0cm 1cm,width=\textwidth]{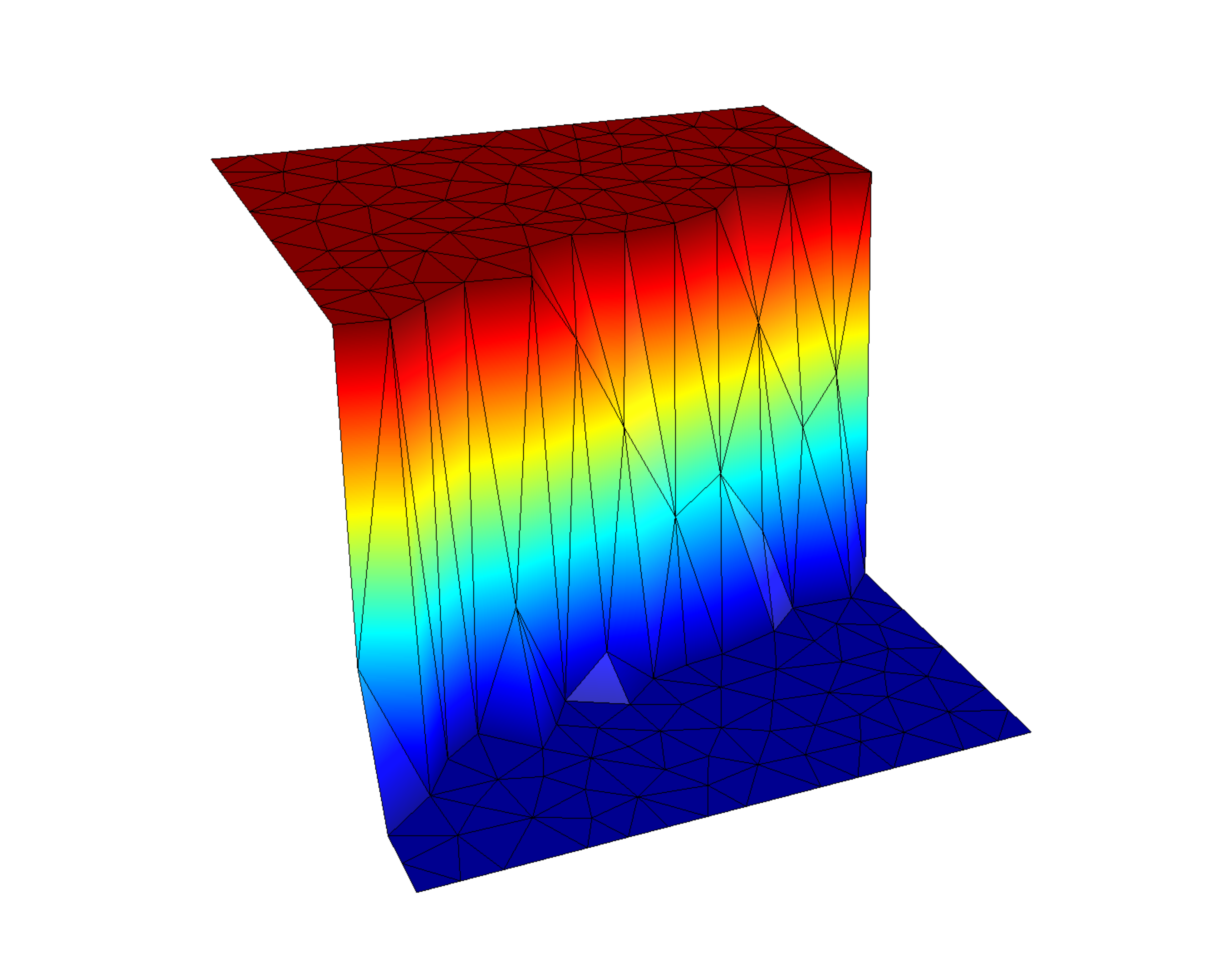}
    \caption{\footnotesize Solution with penalty method.}
    \label{fig:tc1-con}
  \end{subfigure}
  \begin{subfigure}{0.3\textwidth}
  	\vspace{0.2cm}
    \centering
    \includegraphics[width=\textwidth]{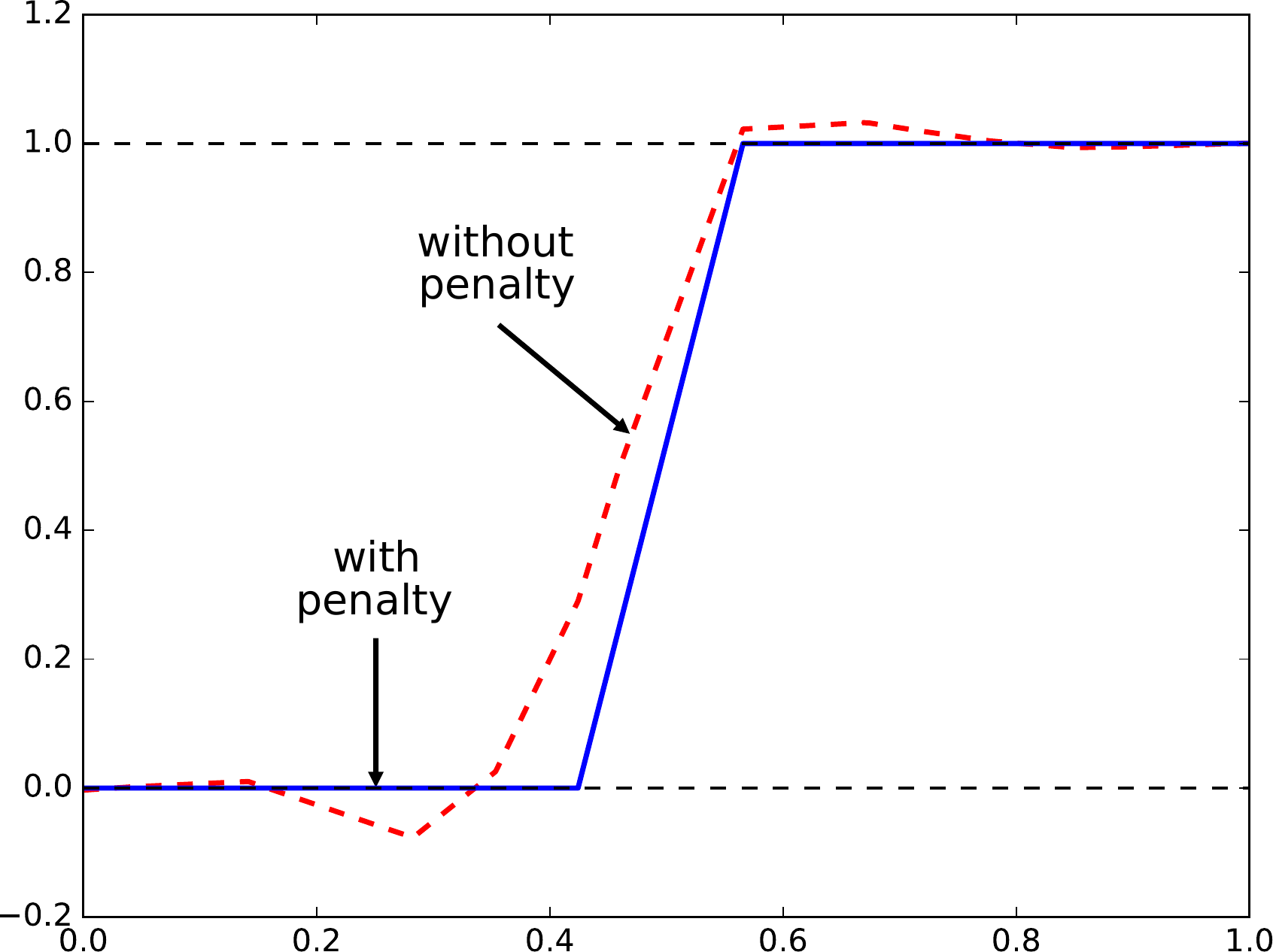}
     \caption{\footnotesize Cross section normal to $\beta$.}
     \label{fig:tc1-cross}
  \end{subfigure}
	\caption{Advection problem over a quasi-uniform mesh.}
	\label{fig:test-case-1_1}
\end{figure}
\begin{figure}[h!]
	\centering
	\begin{subfigure}{0.3\textwidth}
		\centering
		\includegraphics[width=\textwidth]{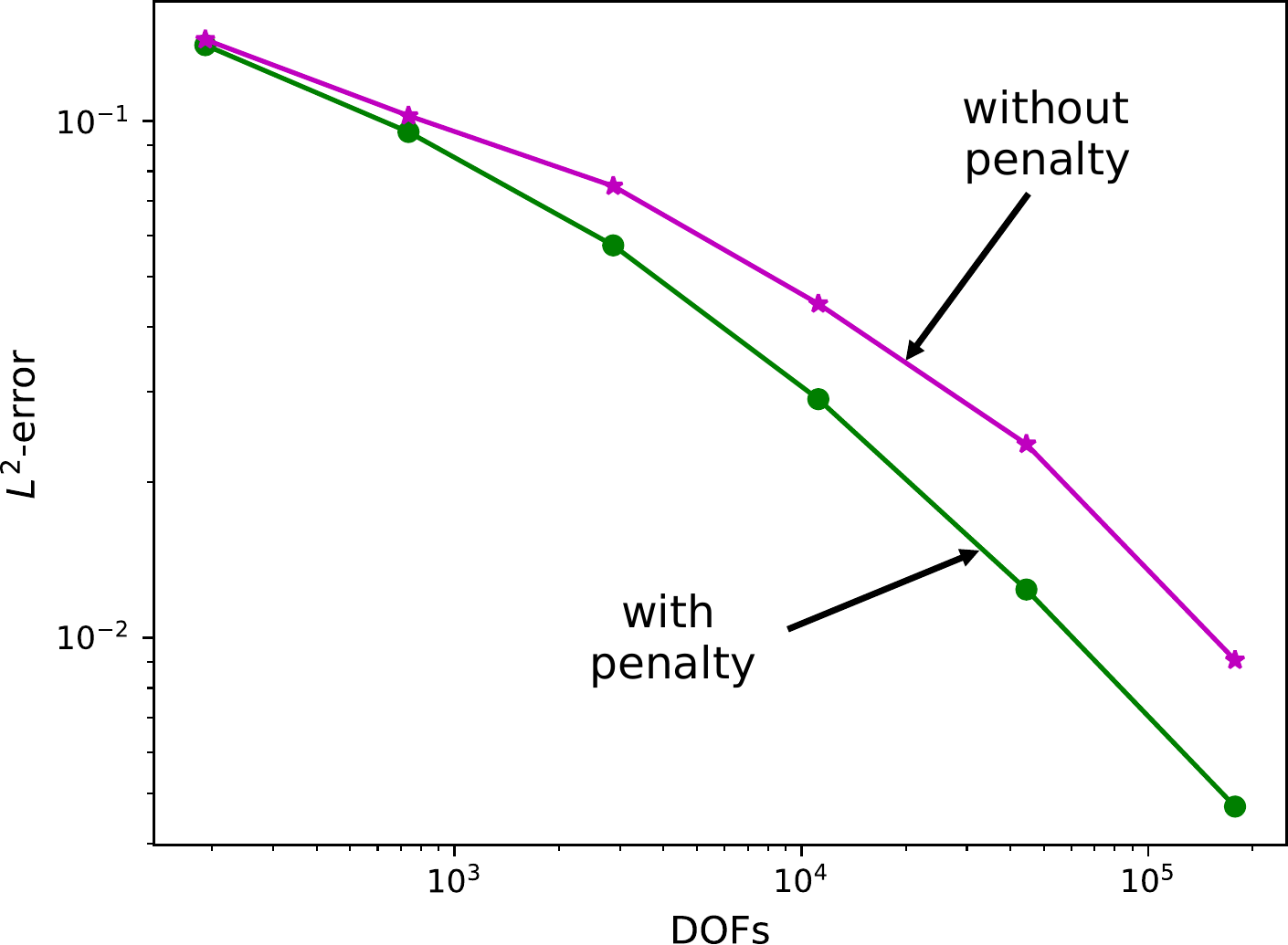}
		\caption{\footnotesize $L^2$ error norm vs DOFs.}
		\label{fig:tc1-l2}
	\end{subfigure}
	\begin{subfigure}{0.3\textwidth}
		\vspace{0.2cm}
		\centering
		\includegraphics[width=\textwidth]{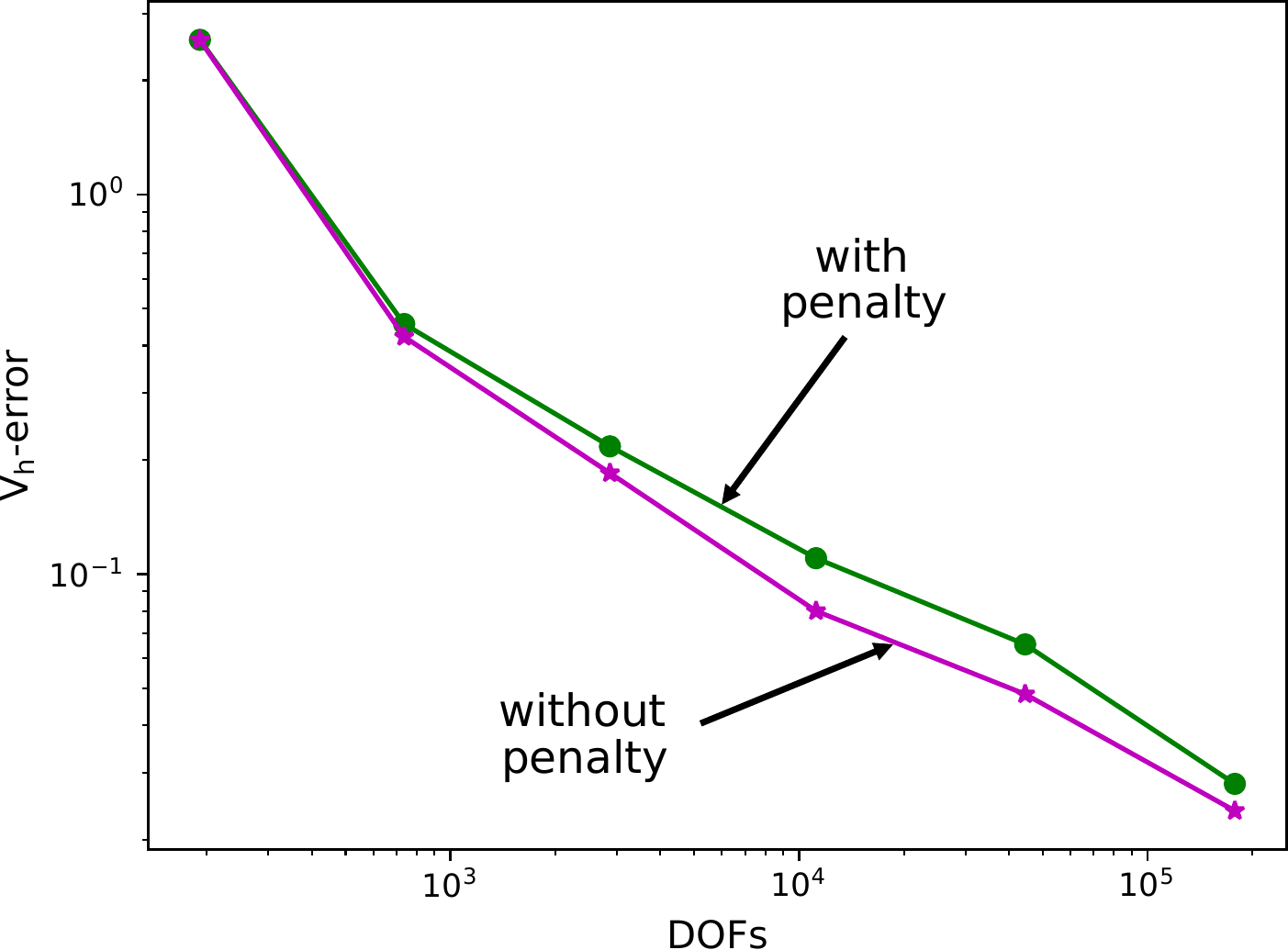}
		\caption{\footnotesize $V_h$ error norm vs DOFs.}
		\label{fig:tc1-vh}
	\end{subfigure}
  \caption{Convergence plots. Uniform refinement.}
  \label{fig:test-case-1_2}
\end{figure}

We simulate a pure advection problem over a quasi-uniform mesh of size ${h = 0.126}$. We set ${\Omega:= (0,1)\times(0,1)}$ and ${\beta=(3/\sqrt{10}, 1/\sqrt{10})^T}$, ${K = 0}$, ${f=0}$. The unit advection field defines that ${\Gamma_-}$ corresponds to the part where ${xy=0}$. The exact solution is ${u=\frac{1}{2}(\tanh((y-\frac{x}{3}-\frac{1}{4})/\epsilon)+1.0)}$, defining an inner layer in the solution of width $\epsilon$. We compute solutions for a sharp layer ($\epsilon=0.01$) using the stabilized method based on residual minimization, both with and without the addition of the nonlinear penalty term. We consider affine ($p=1$) finite elements. Given that the source $f=0$ and the boundary condition $0\leq g \leq 1$ in this experiment, the solution $0 \leq u \leq 1$. Thus, we use the penalty to impose both the lower and upper bounds. Using~\eqref{eq:penalty}, we set $\gamma_0=10^{-5}$. We converge after 18 iterations using ${TOL=10^{-5}}$. As seen in Figures~\ref{fig:tc1-sin}~\&~\ref{fig:tc1-con}, penalties consistently reduce the violation of the solution bounds up to the order of $10^{-3}\%$. Figure~\ref{fig:tc1-cross} shows a cross-section, normal to the advective field. The formulation with penalty significantly improves the bound preservation of the solution, removing the over- and undershoots that appear in the stabilized formulation. Finally, in Figures~\ref{fig:tc1-l2} \&~\ref{fig:tc1-vh}, we show the $L^2$ and $V_h$-error norm convergence, respectively, considering a sequence of uniform meshes. We note that the  constraint enforcement asymptotically produces a worsen convergence in the ${V_h}$-norm, being in line with Remark~\ref{rem:convergence}, while surprisingly exhibiting an improvement in the $L^2$-norm.
\
\subsection{Rotating flow over an adaptive mesh}

\begin{figure}[ht!]
  \centering
  \begin{subfigure}{0.3\textwidth}
    \centering
    \includegraphics[width=\textwidth]{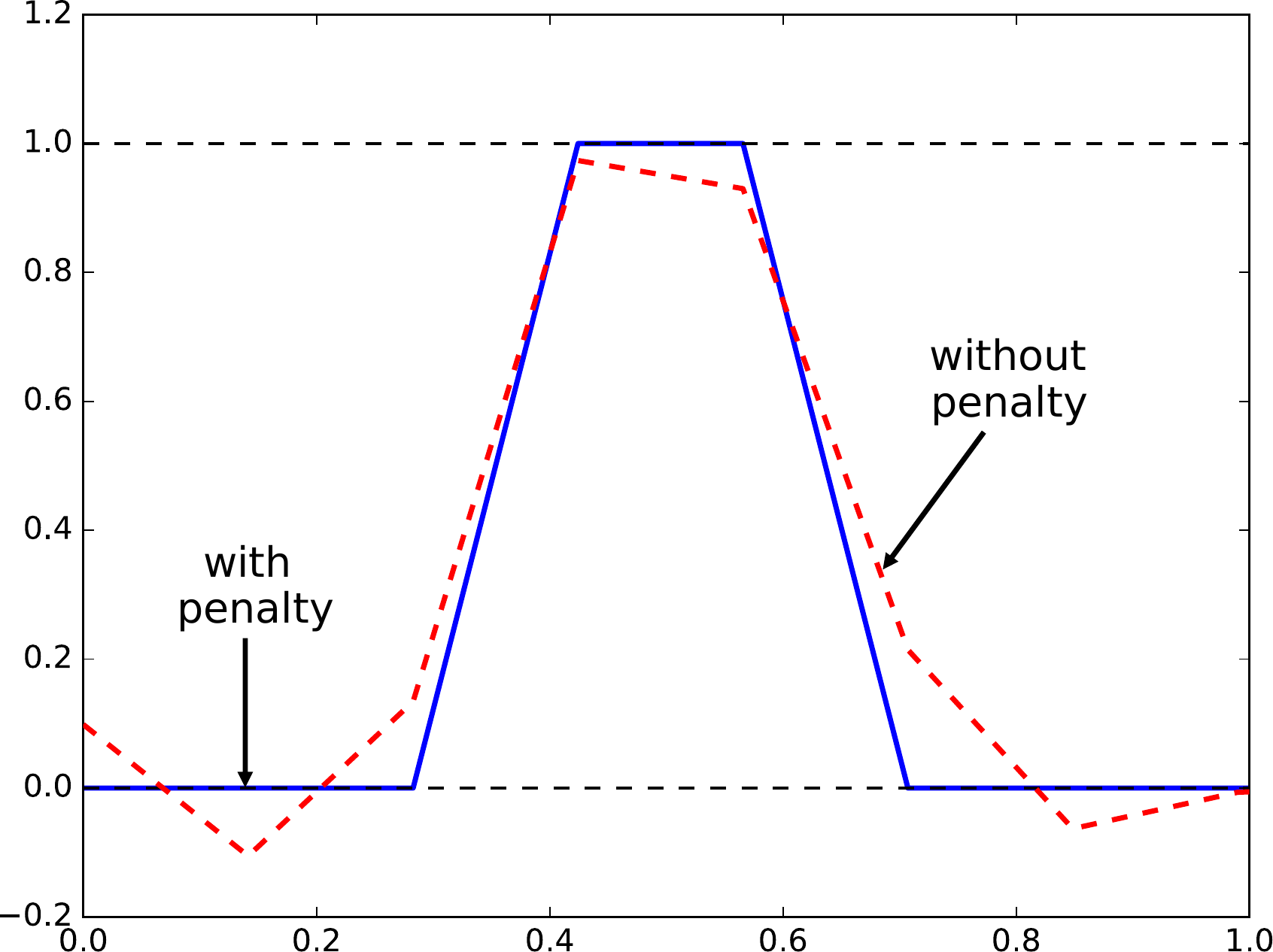}
    \caption{\footnotesize Cross section at $x=y$. Level 10.}
    \label{fig:tc2-cross}
  \end{subfigure}
  \begin{subfigure}{0.33\textwidth}
  	\vspace{0.2cm}
    \centering
    \includegraphics[width=\textwidth]{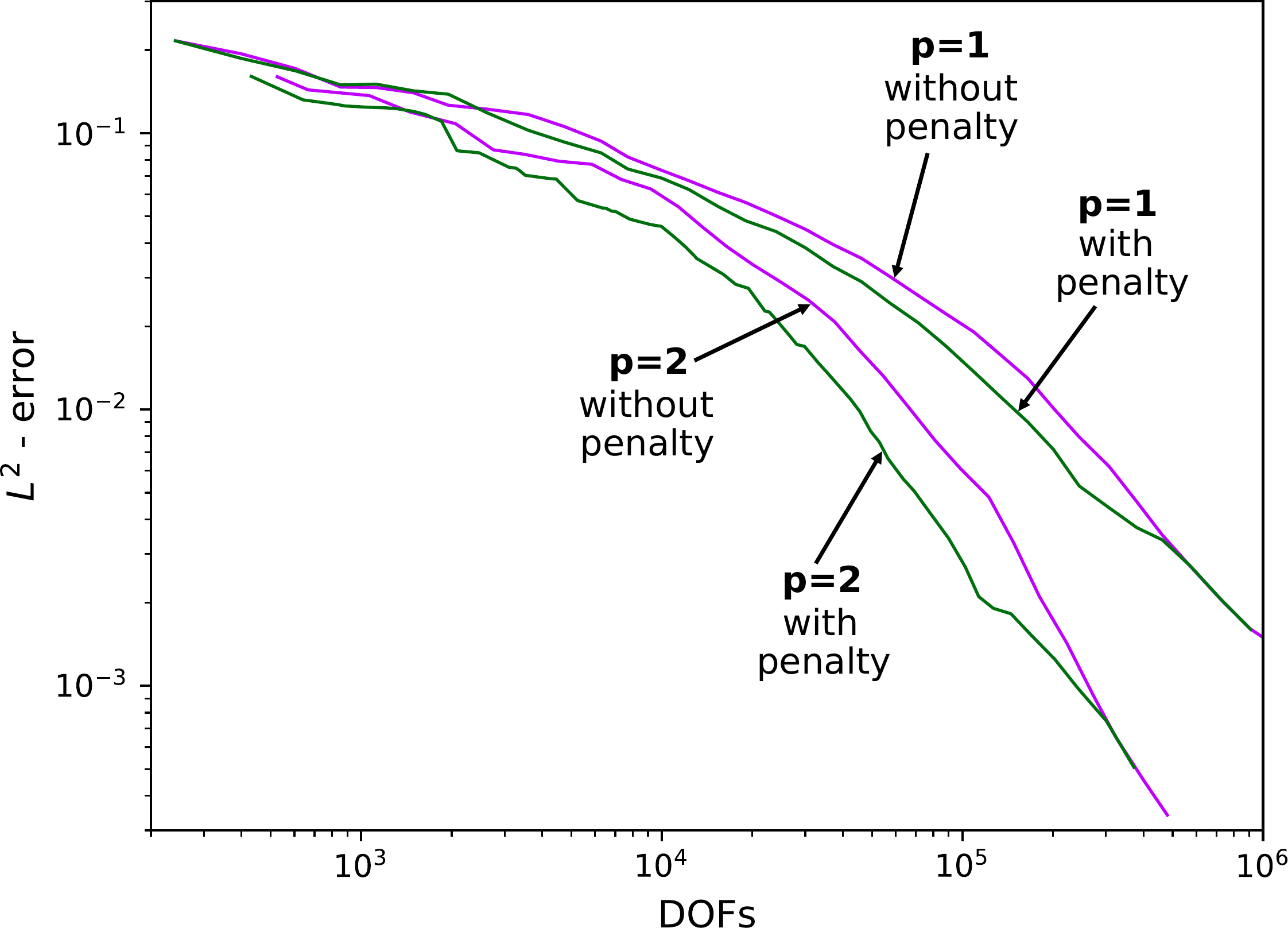}
    \caption{\footnotesize Convergence plot. $L^2$ error norm vs DOFs.}
     \label{fig:tc2-conv}
  \end{subfigure}
  \caption{Rotating flow over an adaptive mesh.}
  \label{fig:test-case-2}
\end{figure}

We now solve a pure-advection test problem proposed in~\cite{kuzmin2010}. Let ${\Omega:= (0,1) \times (-1, 1)}$ with ${\beta=(-y, x)^T}$, ${K = 0}$, ${f=0}$. The convection field rotates counterclockwise, and defines $\Gamma_-=(0,1)\times\{0\} \ \cup \ (0,1)\times \{1\} \ \cup \ \{1\} \times(0,1) \ \cup \ \{0\} \times(-1,0)$. Boundary condition $g$ is:
\begin{equation*}
g =
\left\{
\renewcommand{\arraystretch}{1.4}
\begin{array}{ll}
0.5\{ 1 + \tanh\left[ \epsilon \left(y-0.35\right)\right] \} & \textrm{on} \ (0, 0.5)\times \{0\}, \\
0.5\{ 1 +\tanh\left[\epsilon \left(0.65-y\right)\right] \} & \textrm{on} \ (0.5, 1)\times \{0\}, \\
0 & \textrm{elsewhere on} \ \Gamma_-, \\
\end{array}
\right.
\end{equation*} 
which produces an inner layer in the solution of width $\epsilon$ between 0.35 and 0.65. Similar to the previous test case, we set ${\epsilon=0.01}$. Figure~\ref{fig:tc2-cross} shows a cross-section with and without the inclusion of the penalty term. The bound penalty improves the constraint satisfaction and the inner layer slope. Besides, Figure~\ref{fig:tc2-conv} shows the convergence in $L^2$ and reflects a similar behavior than the uniform mesh case, with the error norm for the penalty formulation solution higher than the one without penalty.

\subsection{Advection-dominated diffusion problem over an adaptive mesh}

\begin{figure}[h!]
  \centering 
  \begin{subfigure}{0.28\textwidth}
    \centering
    \includegraphics[width=\textwidth]{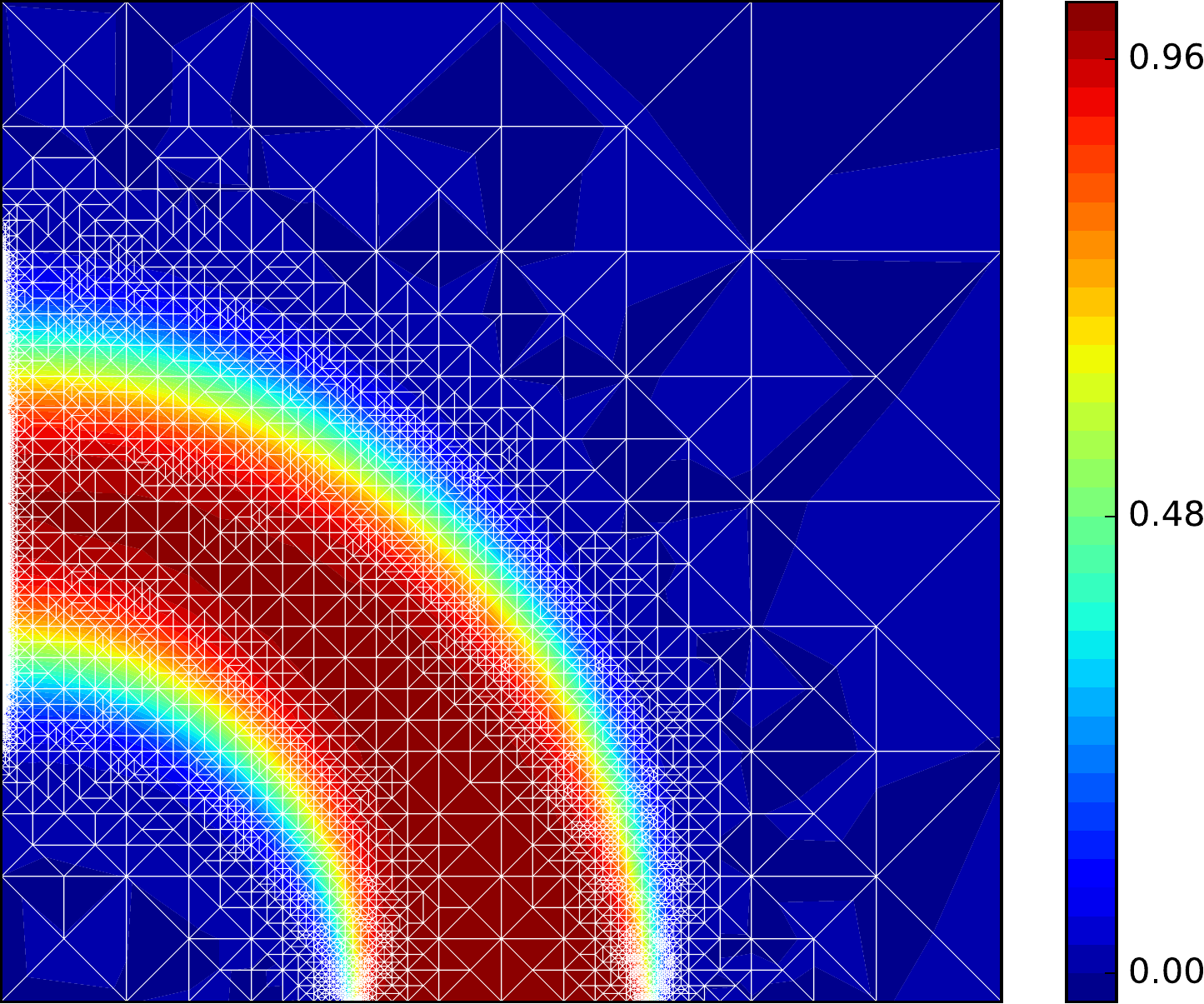} 
    \caption{\footnotesize Computational mesh and 2D solution.}
  \end{subfigure}
  \begin{subfigure}{0.28\textwidth}
  	\vspace{0.2cm}
    \centering
    \includegraphics[width=\textwidth]{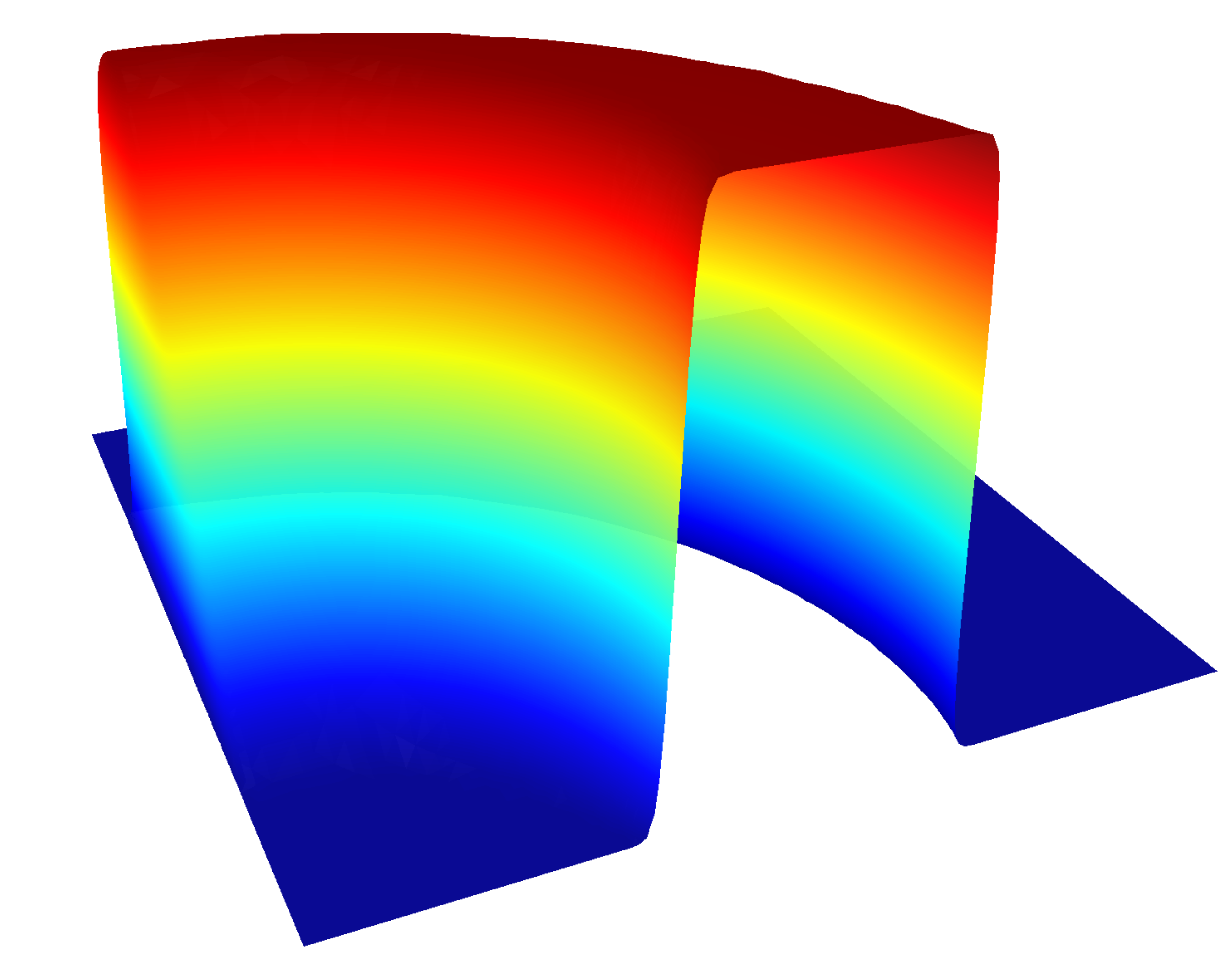}
    \caption{\footnotesize Solution in 3D. Level 25 (82k DOFs).}
  \end{subfigure}
  \caption{Advection-dominated diffusion problem (adaptive mesh).}
  \label{fig:test-case-3}
\end{figure}

We use the nonlinear penalty method to solve a version of the previous test with diffusion. That is, all parameters as above except ${K=10^{-3}}$. This modification induces a boundary layer at ${x=0}$ in the solution due to the contribution of the diffusion part. Our initial mesh is structured and has ${4 \times 4}$ triangular elements. We set ${\gamma_0=10^{-4}}$.  Both trial and test functions are of degree ${p=1}$. The penalty constraints both the lower and upper bounds. Figure~\ref{fig:test-case-3} shows that the adaptive scheme with the nonlinear penalty method captures the boundary layer through a proper error estimate, minimizing the bound violation on each refinement level and thus, delivering physically meaningful solutions at each level.

\section{Conclusions}

We describe a nonlinear weak constraint enforcement for a new adaptive stabilized finite element method. We impose solution bounds on pure-advection and on advection-dominated diffusion problems through the addition of a nonlinear penalty term that weakly enforces the solution range in the variational formulation. The final formulation reduces the bounds violation by several orders of magnitude. Given the stability provided by the formulation, the method moderately increases the computational cost of lower-order schemes.  Finally, this method performs well with adaptive formulations taking advantage of the \textit{a posteriori} error estimate obtained on the fly in the computations. Future work will look for extending the formulation to more complex constraint conditions along with a consistent formulation for transient problems.

\section{Acknowledgments}

This publication was also made possible in part by the CSIRO Professorial Chair in Computational Geoscience at Curtin University and the Deep Earth Imaging Enterprise Future Science Platforms of the Commonwealth Scientific Industrial Research Organisation, CSIRO, of Australia. This project has received funding from the European Union's Horizon 2020 research and innovation programme under the Marie Sklodowska-Curie grant agreement No.~777778 (MATHROCKS). At Curtin University, the Curtin Corrosion Centre, The Institute for Geoscience Research (TIGeR) and by the Curtin Institute for Computation, kindly provide continuing support.  We also acknowledge Alexandre Ern and the anonymous reviewers for the fruitful comments on this work.

\bibliography{ref}

\end{sloppypar}
\end{document}